\documentclass[15pt,a4paper]{article}

\usepackage{amsmath, amsthm, amsfonts, amssymb}
\usepackage{mathptmx}
\usepackage{graphicx}
\usepackage{fancybox}
\usepackage{color}

\numberwithin{equation}{section}
\Large

\def\a{\alpha}

\newtheorem{thm}{Theorem}[section]
\newtheorem{df}{Definition}[section]

\newtheorem{lem}{Lemma}[section]
\newtheorem{cor}{Corollary}[section]
\newtheorem{prop}{Proposition}[section]
\newtheorem{rem}{Remark}[section]
\numberwithin{equation}{section}

\title{One dimensional Weighted Hardy's Inequalities \\and  application. }
\author{Xiaojing Liu,    Toshio Horiuchi, Hiroshi Ando     }
\date{}

\begin{document}
\maketitle

\begin{abstract}

Let   $\Omega$   be  a $C^2$ class bounded domain of $\mathbb R^n\, ( n\ge 1)$.
In the present paper 
we shall improve   one  dimensional  weighted Hardy inequalities with  one-sided boundary condition by adding  sharp remainders.
As  an application, we shall establish  $n$ dimensional  weighted Hardy inequalities with 
  weight functions being  powers of the  distance function $\delta(x)$ to the   boundary $\partial\Omega$.
   Our  results   will  be applicable  to variational problems in a  coming paper \cite{ah2}.
    \footnote{ 
This research was partially supported by Grant-in-Aid for Scientific Research (No. 16K05189) and (No. 15H03621).}

   \end{abstract}
   
   
\section{ INTRODUCTION }

Let $1<p<\infty$  and     $ C_c^\infty((0,1])$ denote the  set of  all $C^\infty$  functions  with   compact  supports in $(0,1]$. 
One  dimensional Hardy inequality with  one-sided boundary condition is  represented by 
\begin{equation} \int_0^1 |u'(t)|^p \,dt\ge \left(1-\frac 1p\right)^p \int_0^1 \frac{|u(t)|^p}{t^p}\,dt+ \left(1-\frac 1p\right)^{p-1} |u(1)|^p \label{A}
\end{equation}
for  every  $u\in C_c^\infty((0,1])$.
When   $u(1)=0$,  this  is a well-known Hardy inequality (see \cite{KO}.  
To  see     the optimality of    coefficient of  the second term in  the  right hand side, by  the  density argument it  suffices  to employ $u_\varepsilon(t)= t^{1-1/p +\varepsilon}$ as  a test  function and make $\varepsilon \downarrow 0$.

Our  first  purpose  in  this  paper is not  only   to establish a weighted version of (\ref{A}) but  also  improve  it  by  adding sharp remainder terms.
As weight functions  we consider power type  weights $t^{\a p}$ for  $t\in [0,1]$.
Surprisingly our   result on  this  matter is  essentially dependent on the  range of  parameter  $\a$.
Let us  explain with symbolic and most simple cases as  examples.
To this  end we classify   the  range  of  the parameter $\a$  into two cases and define the  best constant $\Lambda_{\a,p}$ as  follows:
\begin{df} The  parameter $\a$ is  said  to  be noncritical  and critical if  $\a$  satisfies $\a<1-1/p$  and  $\a\ge  1-1/p$ respectively.
\end{df}
\begin{df}For $1<p<+\infty$ we set
\begin{equation}\Lambda_{\a,p}=\begin{cases} &\left |1-\frac{1}{p}-\a\right |^p,\qquad \text{ if } \a \neq
1-\frac{1}{p},\\   &\left(1-\frac{1}{p}\right)^p,\qquad
\qquad\text{ if } \a = 1-\frac{1}{p}.
\end{cases}\label{2.5}
\end{equation} 
\end{df}
When $\a$  is  noncritical  under  this  definition,  as  a corollary to  Theorem \ref{nct1} we  have a sharp Hardy type  inequality:
\begin{equation} \int_0^1 |u'(t)|^p t^{\a p}\,dt \ge 
\Lambda_{\a,p}\int_0^1 \frac{|u(t)|^p}{t^p}{t^{\a p}}\,dt +
(\Lambda_{\a,p})^{1-1/p} |u(1)|^p,
\end{equation} for  every  $u\in C_c^\infty((0,1])$.
To  see     the optimality of    coefficient of  the second term in  the  right hand side, one  can  employ $u_\varepsilon(t)= t^{1-\a-1/p +\varepsilon}$ as  a test  function as  before.
 When $\a$  is  critical, it  follows  from Proposition \ref{ct1}  that
\begin{equation}\label{C}
\inf_{u\in W}
\int_0^1  |u'(t)|^p t^{\a p} =0,
 \end{equation}
 where $W= \{  u\in C^1([0,1]) :   u(0)=0, u(1)=1\} $. Nevertheless   we  will have  a sharp Hardy type  inequalities (\ref{D}) and  (\ref{E})
 as  a corollary to  Theorem \ref{CT2}.
\par
In Section 2.2, as an important application,  we  will establish  
 $n$ dimensional  weighted Hardy inequalities with 
  weight function being  powers of the  distance function $\delta (x)=\rm{ dist}(x,\partial\Omega)$ to the   boundary $\partial\Omega$. 
  In  this  task it  is  crucial to establish sharp weighted  Hardy  inequalities in  the tubler neighborhood $\Omega_\eta$  of $\Omega$, which are reduced  to  the  one  dimensional inequalities  in Section 2.1.
To  this  end   $\Omega $ is  assumed  to be  a bounded domain of  $\mathbf R^N$ ( $N\ge 1$ ) whose boundary  $\partial\Omega$ is a $C^2$ compact manifolds in  the  present  paper. 
We  prepare  more  notations to  describe  our  results.
For   $ \alpha \in \bf {R}$,
by
  $ L^p(\Omega,\delta^{p\a})$ we denote the space of Lebesgue measurable functions with weight $ \delta^{\a p}$, 
for which 
\begin{equation} || u ||_{ L^p(\Omega,\delta^{p\a})} = \bigg( \int_{\Omega}|u|^p \delta^{\a p}\, 
dx\bigg ) ^{1/p} < +\infty.\label{2.1} \end{equation} 
 $W_{\a,0}^{1,p}(\Omega)$ is given by the  completion of $C_c^\infty(\Omega)$ with respect to the norm 
defined by
\begin{equation} 
|| u||_{ W_{\a,0}^{1,p}(\Omega)} =
|| |\nabla u|  ||_{ L^p(\Omega, \delta ^{p\a})} + || u||_{ L^p(\Omega,\delta^{p\a})}. \label{2.2}
\end{equation}
Then 
 $W^{1,p}_{\a,0}(\Omega) $  becomes a Banach space  with the norm $|| \cdot||_{ W_{\a,0}^{1,p}(\Omega)}
$.
 Under these preparation 
we will  state  the  noncritical weighted Hardy inequality as
Theorem \ref{NCT1}, which is    the   counter-part  to  Theorem \ref{nct1}.
In particular as its corollary,  we   have the simplest one:
\begin{equation}\label{HI}
\int_\Omega  |\nabla u|^p \delta^{\alpha p} \ge \mu \int_\Omega |u|^p \delta ^{ p(\alpha -1)}, \qquad \forall u\in W^{1,p}_{\alpha,0}(\Omega),
 \end{equation}
where   $\a< 1-\frac 1p$  and $\mu$ is a positive  constant   essentially depending on the boundary $\partial\Omega$. 
If  $\a=0$ and  $p=2$, then  (\ref{HI})  is  a well-known Hardy inequality and  valid for  a bounded domain $\Omega $  of  $\mathbf R^N$ with Lipschitz boundary  (c.f. \cite{BM,D, MMP}).
Further
if     $\Omega$ is  convex  and  $ \a=0$,  then   $\mu = \Lambda_{0,p}$ holds  for   arbitrary $1<p<\infty$ (see \cite{MS}). \par 
It  is  worthy to remark that (\ref{HI}) is  never valid in  the critical    case that $\a\ge  1-1/p$
by (\ref{C}) ( see also  Proposition \ref{CT1} ). Nevertheless,
 we  will   establish   in  this   case a  variant of 
 weighted  Hardy's inequalities as Theorem \ref{CT3} which correspond  to those  in  Theorem \ref{CT2}.
 As  its  corollary we describe     Hardy's  inequalities  with  a compact perturbation which  are  closely relating   to  
 the so-called weak Hardy property of  $\Omega$.  We  remark  that
 a constant $\gamma^{-1}$ in (\ref{2.14})  and (\ref{2.15}) concerns  the weak Hardy constant,  but  in  this  case   the  strong Hardy  constant  is $+\infty$ (  see \cite{D} for  the  detail).
  In \cite{ah0},   two  of  the  authors   have improved   the  weighted Hardy inequalities
   adopting  $|x|^{\a p}$  ( 
powers  of  distance to the  origin $O\in \Omega$ )   as   weight  functions   instead of  $\delta^{\a p}$.
In  the  present paper,  some  inequalities 
of  Hardy type  in \cite{ah0} are   employed with minor  modifications, especially  when $1<p<2$
 (see also \cite{ah1,dha,dha2,aha3}).
 We  note  that our results will be further improved in \cite{H} for
 non-doubling weights.
Lastly we    remark  that our  results   will  be applicable  to variational problems in a  coming paper \cite{ah2}.

\par
 This  paper is  organized in  the  following way: The main  results are described  in Section 2. 
Theorem 2.1 and Theorem 2.2 are  established in Section 3.  Theorem 2.3 and Theorem 2.5 together with   their corollaries are proved in Section 4. 
The  proof  of Theorem 2.4 is given in Section 5 and the  proofs of  Proposition  2.1 and Proposition 2.2 are  given  in Section 6.
In  Appendix  the  proofs of Lemma 3.2 and  Lemma 3.4 are provided for  the  sake of  self-containedness.
 
\section{Main results  }
\begin{df} \label{def1.1}For $t\in (0,1)$  and $R> e$,  we  set
\begin{equation}
A_1(t):=\log\frac{R}{t},\qquad
A_2(t):=\log A_{1}(t)
.\label{1.3}
\end{equation}\end{df}
\subsection{Results in the  one  dimensional  case}
The  proofs of Theorem \ref{nct1} and Theorem \ref{CT2} including corollaries will be  given in Section 3 and Appendix.
\begin{thm}\mbox{( \rm Noncritical case )} \label{nct1}
Assume that $ \a<1-1/p$, $1<p<\infty$ and  $R>e$. Then, there exist  positive numbers  $C_0=C_0(\a,p,R )$,  $C_1=C_1(\a,p,R )$ and  $L=L(\a,p,R )$ such  that  for  every $u\in C_c^\infty((0,1])$, we  have 
\begin{align}  \int_0^1  & \left( |u'|^p -\left|\frac{u}{t}\right|^p \left(\Lambda_{\a,p}    + \frac{C_0 }{A_1(t)^2}\right)\right) t^{\a p}\,dt  \notag \\ &\ge  C_1
 \int_0^1 \left( |u'|^p + \left|\frac{u}{t}\right|^p \left(\Lambda_{\a,p} + \frac{C_0 }{A_1(t)^2}\right) \right) t^{\a p+1}\,dt +L|u(1)|^p.\label{nc1}
\end{align}
\end{thm}
\begin{cor}\label{cor2.1}
Assume that $ \a<1-1/p$ and $1<p<\infty$. Then, for  every  $u\in C_c^\infty((0,1])$
\begin{equation} \int_0^1 |u'(t)|^p t^{\a p}\,dt \ge 
\Lambda_{\a,p}\int_0^1 \frac{|u(t)|^p}{t^p}{t^{\a p}}\,dt +
({\Lambda_{\a,p}})^{1-1/p} |u(1)|^p.\label{B}
\end{equation} 
\end{cor}
In  the  critical case  we  have   somewhat more precise results.
\begin{thm}\mbox{\rm ( Critical case)} \label{CT2}
\begin{enumerate} \item
Assume that $ \a>1-1/p$, $1<p<\infty$ and $R>e$. Then there exist  positive numbers    $C_0=C_0(\a,p,R )$,  $C_1=C_1(\a,p,R )$ and  $L=L(\a,p,R )$  such  that  for  every $u\in C_c^\infty((0,1])$, we  have 
\begin{align}  \int_0^1  & \left( |u'|^p -\left|\frac{u}{t}\right|^p \left(\Lambda_{\a,p}    + \frac{C_0 }{A_1(t)^2}\right)\right) t^{\a p}\,dt  +L |u(1)|^p \notag \\ &\ge  C_1
 \int_0^1 \left( |u'|^p + \left|\frac{u}{t}\right|^p \left(\Lambda_{\a,p} + \frac{C_0 }{A_1(t)^2}\right) \right) t^{\a p+1}\,dt.\label{c1}
\end{align}
\item 
Assume that $ \a=1-1/p$, $1<p<\infty$ and $R>e^e$. Then, there exist  positive numbers $C_0=C_0(\a,p,R )$,  $C_1=C_1(\a,p,R )$and  and  $L=L(\a,p,R )$ such  that  for  every $u\in C_c^\infty((0,1])$, we  have 
\begin{align}  \int_0^1  & \left( |u'|^p -  \left|\frac{u}{t}\right|^p\frac{1}{A_1(t)^p} \left(\Lambda_{\a,p}   + \frac{C_0 }{A_2(t)^2}\right)\right) t^{p-1}\,dt +L |u(1)|^p  \notag \\ &\ge  C_1
 \int_0^1 \left( |u'|^p + \left|\frac{u}{t}\right|^p\frac{1}{A_1(t)^p}  \left(\Lambda_{\a,p}   + \frac{C_0 }{A_2(t)^2}\right) \right) t^{p}\,dt. \label{2.10}
\end{align}
\end{enumerate}
\end{thm}
\begin{cor}\label{cor2.2}
 \begin{enumerate}\item If  $\a>1-1/p$ and $1<p<\infty$,  then for  every  $u\in C_c^\infty((0,1])$
\begin{equation} \int_0^1 |u'(t)|^p t^{\a p}\,dt + (\Lambda_{\a,p})^{1-1/p} |u(1)|^p\ge\Lambda_{\a,p} \int_0^1 \frac{|u(t)|^p}{t^p}{t^{\a p}}\,dt. \label{D}
\end{equation}
 \item 
If $\a=1-1/p$, $1<p<\infty$ and $R>e$,  then  for  every  $u\in C_c^\infty((0,1])$
\begin{equation} \int_0^1 |u'(t)|^p t^{p-1}\,dt +(\Lambda_{\a,p})^{\a} A_1(1)^{1-p} |u(1)|^p\ge\Lambda_{\a,p} \int_0^1 \frac{|u(t)|^p}{t A_1(t)^p}\,dt. \label{E}
\end{equation} 
\end{enumerate}
\end{cor}
\begin{rem}
We remark  that   Corollaries \ref{cor2.1} and \ref{cor2.2} follow direct from the  arguments in  the  proofs of the  corresponding theorems 
except  for the optimality of  the  constant $L =( \Lambda _{\a,p})^{1-1/p}$. For  the  proofs  of the optimality, one  can employ as test functions $u_{\varepsilon}= t^{1-\a-1/p+\varepsilon} $ in (\ref{B}), $u_{\varepsilon}= t^{1-\a-1/p+\varepsilon} $ in (\ref{D})  and $u_{\varepsilon}=
 A_1(t)^{1-1/p-\varepsilon} $ in (\ref{E})  respectively  with $\varepsilon$  being   sufficiently small.
\end{rem}
Further we  remark  an  elementary result which will  be  useful  in  the  subsequent.
\begin{prop}\label{ct1} $( $\rm {Critical  case} $)$
Assume that $ \a\ge 1-1/p$ and $1<p<\infty$. Then  we have 
\begin{equation}\label{}
\inf_{u\in W}
\int_0^1  |u'(t)|^p t^{\a p} =0,
 \end{equation}
 where $W= \{  u\in C^1([0,1]) :   u(0)=0, u(1)=1\}.  $
 \end{prop}
The proof  will be  given in Section 6.

\subsection{Results in a  domain  of  $\mathbf R^N $ }
The  proofs of  Theorem \ref{NCT1}, Corollary \ref{NCC1}, and Theorem \ref{CT3} will  be given in Section 4. Theorem \ref{NCT2} will be  proved in Section 5.
 Let $\delta(x)=  \rm{dist}(x,\partial\Omega)$. We  use the  following notations:
 \begin{equation}
\Omega_\eta = \{ x\in \Omega:  \delta(x)<\eta \}, \qquad
\Sigma_\eta = \{ x\in \Omega:  \delta(x)=\eta \}. \label{NBD}
\end{equation}

\begin{thm}\mbox{ \rm ( Noncritical case )} \label{NCT1} Assume  that  $\Omega$  is  a  bounded domain of  class $C^2$ in $\mathbf R^N$.
Assume  that   $ \a<1-1/p$, $1<p<\infty$  and $R> e\cdot\sup_{x\in \Omega}\delta(x)$. Assume  that  $\eta $  is  a sufficienty  small positive  number. Then, there exist  positive  numbers  $C_2=C_2(\a,p,R,\eta)$ and $L=L(\a,p,R,\eta)$  such  that  for every  $ u\in W^{1,p}_{\alpha,0}(\Omega)$, we have 
\begin{equation}\label{2.6}
\int_{\Omega_\eta} \left( |\nabla u|^p  -\Lambda_{\a,p}\left| \frac{u}{\delta} \right|^p \right)\delta^{\alpha p} 
\ge C_2 \int_{\Omega_\eta} \left | \frac{u}{\delta}\right |^p \frac{1}{ A_1(\delta)^2}\delta ^{p\alpha } +L \int_{\Sigma_\eta} |u|^p\delta^{\a p}.
 \end{equation}
 \end{thm}

\begin{cor}\label{NCC1}
Under the  same assumptions  as in Theorem \ref{NCT1}, there exists  a positive number  $\gamma=\gamma(\a,p,R)$   such  that for every  $ u\in W^{1,p}_{\alpha,0}(\Omega)$, we have 
\begin{equation}\label{2.7}
\int_{\Omega} \left( |\nabla u|^p  -\gamma\left| \frac{u}{\delta} \right|^p \right)\delta^{\alpha p} 
\ge 0. 
 \end{equation}
 
 \end{cor}
Moreover 
for  any  bounded  domain $\Omega \subset \mathbf R^N$
 we can prove  the following:
 \begin{thm}\label{NCT2}\mbox{ \rm(  Noncritical case )}  Assume  that  $\Omega$  is  a  bounded domain of   class $C^2$ in $\mathbf R^N$.
Assume  that $ \a<1-1/p$, $1<p<\infty$  and $R> e\cdot\sup_{x\in \Omega}\delta(x)$.  Then, the followings are equivalent.
\begin{enumerate}
\item  There exists  positive  a number  $\gamma$   such  that the  inequality (\ref{2.7}) is valid for  every  $u\in W^{1,p}_{\alpha,0}(\Omega)$.

\item  For  a sufficiently  small  $\eta>0$,
there exist  positive  numbers  $\kappa$, $C_2$ and $L$ such  that the  inequality (\ref{2.6})  with   $\Lambda_{\a,p}$  replaced by  $\kappa$ is  valid for   every  $u\in W^{1,p}_{\alpha,0}(\Omega)$.
\end{enumerate} 
\end{thm}

\par\medskip
\begin{thm}\label{CT3} \mbox{\rm ( Critical case)} Assume  that  $\Omega$  is  a  bounded domain of  class $C^2$ in $\mathbf R^N$.
\begin{enumerate}
\item
Assume that $ \a>1-1/p$, $1<p<\infty$ and $R> e^e\cdot\sup_{x\in \Omega}\delta(x)$. 
Assume  that  $\eta $  is  a sufficienty  small positive  number. 
Then, there exist  positive  numbers    $C_2=C_2(\a,p,R,\eta)$   and  $L=L(\a,p,R,\eta)$ such  that  
for every  $ u\in W^{1,p}_{\alpha,0}(\Omega)$, we have 
\begin{align}\label{2.11}
\int_{\Omega_\eta}   \left( |\nabla u|^p  -\Lambda_{\a,p}\left| \frac{u}{\delta} \right|^p \right)\delta^{\alpha p} + L \int_{\Sigma_\eta} |u|^p\delta^{\a p}
\ge C_2 \int_{\Omega_\eta} \left | \frac{u}{\delta}\right |^p \frac{\delta ^{p\alpha }}{ A_1(\delta)^2}.
 \end{align}
\item
Assume  that  $ \a=1-1/p$,$1<p<\infty$ and $R> e^e\cdot\sup_{x\in \Omega}\delta(x)$.  
Assume  that  $\eta $  is  a sufficienty  small positive  number. 
Then, there exist  positive  numbers  $C_2=C_2(\a,p,R,\eta)$  and  $L=L(\a,p,R,\eta)$
 such  that  for every  $ u\in W^{1,p}_{\alpha,0}(\Omega)$, we have 
\begin{align}\label{2.12}
\int_{\Omega_\eta} & \left( |\nabla u|^p  -\Lambda_{\a,p}\left| \frac{u}{\delta}  \right|^p  \frac{1}{A_1(\delta)^p}\right)\delta^{p-1} + L \int_{\Sigma_\eta} |u|^p\delta^{p-1} \\ &
\ge C_2 \int_{\Omega_\eta} \left | \frac{u}{\delta}\right |^p \frac{1}{ A_1(\delta)^p}\frac{1}{A_2(\delta)^2} \delta ^{p-1}.
\notag  \end{align}
\end{enumerate}
\end{thm}

\begin{cor} \label{CCT3}Assume  that  $\Omega$  is  a  bounded domain of  class $C^2$ in $\mathbf R^N$.
\begin{enumerate}
\item
Assume that $ \a>1-1/p$, $1<p<\infty$  and $\gamma=\gamma(\a,p,R,\eta)$. 
Then, there exists a  positive  number   $L'=L'(\a,p,R,\eta)$ such  that  for every  $ u\in W^{1,p}_{\alpha,0}(\Omega)$, we have
\begin{align}\label{2.14}
\int_{\Omega}  & \left( |\nabla u|^p  -\gamma \left| \frac{u}{\delta} \right|^p \right)\delta^{\alpha p} + L' \int_{\Sigma_\eta} |u|^p\delta^{\a p}
\ge0.
 \end{align}
\item
Assume  that  $ \a=1-1/p$, $1<p<\infty$ and $R> e^e\cdot\sup_{x\in \Omega}\delta(x)$.  
 Then, there exists  positive  numbers  $\gamma$  and  $L'=L'(\a,p,R,\eta)$
 such  that  for every  $ u\in W^{1,p}_{\alpha,0}(\Omega)$, we have
\begin{align}\label{2.15}
\int_{\Omega} \left( |\nabla u|^p  -\gamma\left| \frac{u}{\delta}  \right|^p  \frac{1}{A_1(\delta)^p}\right)\delta^{p-1} + L' \int_{\Sigma_\eta} |u|^p\delta^{p-1} 
\ge 0.
  \end{align}
\end{enumerate}
\end{cor}

\begin{prop}\label{CT1} $( $\rm {Critical  case} $)$ Assume  that  $\Omega$  is  a  bounded domain of  class $C^2$ in $\mathbf R^N$.
Assume  that   $1<p<\infty$  and  $ \a \ge1-1/p$. 
Then for  arbitrary  $\eta \in (0, \sup_{x\in \Omega} \delta(x) )$ we have 
\begin{equation}\label{2.8}
\inf
\left\{
\int_{\Omega}  |\nabla u|^p \delta^{\a p} : u\in C^1_c (\Omega),  u=1 \mbox{  on } \{ \delta(x)=\eta\}\right \}=0.
 \end{equation}
 \end{prop}
 The proof  will be  given in Section 6.
From this
it  is worthy to remark that Hardy's inequality (\ref{HI}) never  holds  in  the  critical case.

\section{Proofs of Theorem \ref{nct1}  and Theorem  \ref{CT2} }

\subsection{Auxiliary inequalities  in  the noncritical case}
When $p=2$ and $\a=0$, the first  lemma is  established in \cite{BM}. 
\begin{lem}\label{l2} \mbox{\rm ( Noncritical case)}
Assume  that  $ f\in C^{}([0,1])\cap C^1((0,1])$ is a monotone nondecreasing function such  that $f(1)\le 1$.  Assume that $ \a<1-1/p$ and $1<p<\infty$.
Then for  every
  $u\in C_c^\infty((0,1])$, we  have 
\begin{equation} \int_0^1\left ( |u'|^p -\Lambda_{\a,p} \left|\frac{u}{t}\right|^p \right)t^{\a p}\,dt \ge  
\int_0^1 \left( |u'|^p -\Lambda_{\a,p}\left|\frac{u}{t}\right|^p \right)t^{\a p} f  \,dt. \label{3.2}
\end{equation}
In particular we  have 
\begin{equation} \int_0^1 |u'|^p t^{\a p}\,dt \ge  
\Lambda_{\a,p} \int_0^1\left|\frac{u}{t}\right|^p t^{\a p}  \,dt. 
\end{equation}

\end{lem}
 \par\medskip
\noindent{\bf Proof of Lemma \ref{l2}:} Without loss of generality 
we  assume  that $f\ge 0$,  $f(1)=1$, and $u\ge0$.  Define  $ g= 1-f$. Then  $g\ge 0$  and  $ g'\le 0$.
By integration by parts we have 
\begin{align*} (1-\a-1/p) & \int_0^1 u^p t^{p(\a-1)} g\,dt = -1/p \left[ u^p t^{p(\a-1)+1} g\right]_0^1 \\ &
+ 1/p \int_0^1 u^p t^{p(\a-1)+1} g'\,dt + \int_0^1 u^{p-1}  u' t^{p(\a-1)+1} g\,dt.
\end{align*}
Since $g'=-f'\le 0$  and  $g\ge 0$,
\begin{align*} (1-\a-1/p) \int_0^1 u^p t^{p(\a-1)} g\,dt\le \int_0^1 u^{p-1}  u'  t^{p(\a-1)+1} g\,dt.
\end{align*}
By  H\"older's inequality, noting  that ${p(\a-1)+1}= {(p-1)(\a-1)+\a}$,  we have 

\begin{align*} (1-\a-1/p)\left( \int_0^1 u^p t^{p(\a-1)} g\,dt \right)^{1/p}\le \left(\int_0^1  |u' |^pt^{p\a} g\,dt\right)^{1/p}.
\end{align*}
Using $ g=1-f$ and  the definition of  $\Lambda_{\a,p}$, we have 
\begin{equation*} \int_0^1 ( |u'|^p -\Lambda_{\a,p}({u}/{t})^p)t^{\a p}\,dt \ge  
\int_0^1 ( |u'|^p -\Lambda_{\a,p}({u}/{t})^p)t^{\a p} f  \,dt. 
\end{equation*}
\hfill $\Box$
\begin{lem}\label{l3}  \mbox{\rm ( Noncritical case)} Assume that $ \a<1-1/p$,$1<p<\infty$ and $R>e$. Then, there exist  positive numbers $C_3=C_3(\a,p,R )$ and  $L=L(\a,p,R )$ such  that  for  every $u\in C_c^\infty((0,1])$, we  have 
\begin{equation} \int_0^1 \left ( |u'|^p -\Lambda_{\a,p} \left|\frac{u}{t}\right|^p \right)t^{\a p}\,dt \ge  C_3
\int_0^1  \left|\frac{u}{t}\right|^p t^{\a p} \frac{1}{A_1(t)^{2} } \,dt  + L|u(1)|^p.\label{3.3}
\end{equation}
\end{lem}

The  proofs  of  Lemma \ref{l3} together with  Lemma \ref{l7}   will be   given in \S 6.
It  follows  from  Lenma \ref{l2} and  Lemma \ref{l3} that  we  have 
\begin{lem}\label{l4}  \mbox{\rm ( Noncritical case)} Assume that $ \a<1-1/p$, $1<p<\infty$ and $R>e$. There exist  positive numbers $C_4=C_4(\a,p,R )$ and  $L=L(\a,p,R )$ such  that  for  every $u\in C_c^\infty((0,1])$, we  have 
\begin{equation} \int_0^1\left ( |u'|^p -\Lambda_{\a,p} \left|\frac{u}{t}\right|^p \right )t^{\a p}\,dt \ge  C_4
\int_0^1\left ( |u'|^p +\Lambda_{\a,p} \left|\frac{u}{t}\right|^p \right) \frac{t^{\a p} }{A_1^{2} } \,dt  + L|u(1)|^p.\label{3.4}
\end{equation}
In  particular $C_4$ is  given  by
$ C_4= C_3/(1+2\Lambda_{\a,p} ).$
\end{lem}

\noindent{\bf Proof of Lemma \ref{l4}}
We use   Lemma \ref{l2} for  $f= C_3 A_1^{-2} $ with   $C_3$ being small.  Multiplying (\ref{3.3}) by $2\Lambda_{\a,p} $ and adding it  to (\ref{3.2}), we have (\ref{3.4}) for  
$C_4= C_3/(1+2\Lambda_{\a,p} )$. \par \hfill $\Box$

\subsection{Proof of Theorem \ref{nct1}}
By  adding $$-C_0\int_0^1 \left|\frac{u}{t}\right|^p \frac{1}{A_1(t)^2} t^{\a p}\,dt$$ to  the  both side
of  
(\ref{3.4}) we  have  
\begin{align}  \int_0^1  & \left( |u'|^p -\left|\frac{u}{t}\right|^p \left(\Lambda_{\a,p}    + \frac{C_0 }{A_1(t)^2}\right) \right) t^{\a p}\,dt  \\ 
&\ge  C_4
 \int_0^1 \left( |u'|^p + \left|\frac{u}{t}\right|^p  \left( \Lambda_{\a,p}-\frac{C_0}{C_4} \right)\right) \frac{1}{A_1(t)^2} t^{\a p}\,dt  +L|u(1)|^p.\notag
\end{align}
Now  we set $\,  C_0= {\Lambda_{\a,p} C_4}/{3} ,\, C_1'={ C_4}/{3}.$
Assuming  that  $C_0\le \Lambda_{\a,p} (\log R)^2$,  we  have 
$ C_0/A_1(t)^2 \le  C_0 /( \log R)^2\le \Lambda_{\a,p}$,  and  hence  
$$ \Lambda_{\a,p}-\frac{C_0}{C_4} = \frac23 \Lambda_{\a,p}\ge \frac 13 \left(\Lambda_{\a,p}  + \frac{C_0}{A_1(t)^2}\right).$$
 Then we  have 
\begin{align}  \int_0^1  & \left( |u'|^p -\left|\frac{u}{t}\right|^p \left(\Lambda_{\a,p}    + \frac{C_0 }{A_1(t)^2}\right) \right) t^{\a p}\,dt \notag \\ 
&\ge  C_1'
 \int_0^1 \left( |u'|^p + \left|\frac{u}{t}\right|^p  \left( \Lambda_{\a,p}+\frac{C_0}{A_1(t)^2} \right)\right) \frac{1}{A_1(t)^2} t^{\a p}\,dt +L|u(1)|^p.\notag
\end{align}
By a  calculation  we see that    $  t A_1^{2} \le 4R/e^2  \,\, ( t\in [0,1]).$  
Thus  the  desired  inequality holds  for  $C_1= C_1' e^2/(4R)$.
 \hfill $\Box$
 
\subsection{Auxiliary inequalities in  the critical case}
The following lemma will  be  established in Section 6 together  with Lemma \ref{l3}.
\begin{lem}\mbox{\rm ( Critical case)} \label{l7} 
\begin{enumerate}
\item
 Assume that $\a>1-1/p$, $1<p<\infty$ and $R>e$. There exist  positive numbers  $C_5=C_5(\a,p,R )$ and  $L=L(\a,p,R )$ such  that  for  every $u\in C_c^\infty((0,1])$, we  have
 \begin{align} \int_0^1  \left ( |u'|^p -\Lambda_{\a,p} \left|\frac{u}{t}\right|^p  \right) & t^{\a p}\,dt +L |u(1)|^p 
\ge  C_5
\int_0^1  \left|\frac{u}{t}\right|^p \frac{ t^{\a p}}{A_1^2} \,dt. \label{3.7} 
\end{align}
\item
 Assume that $\a=1-1/p$, $1<p<\infty$ and $R>e$. There exist  positive numbers $C_5=C_5(\a,p,R )$ and  $L=L(\a,p,R )$ such  that  for  every $u\in C_c^\infty((0,1])$, we  have
\begin{align} \int_0^1  \left ( |u'|^p -\Lambda_{\a,p} \left|\frac{u}{t}\right|^p \frac{1}{A_1^p} \right)  t^{p-1}\,dt +L |u(1)|^p
\ge  C_5
\int_0^1  \left|\frac{u}{t}\right|^p \frac{ t^{p-1}  }{A_1^p A_2^2} \,dt.  \label{3.8} 
 \end{align}
\end{enumerate}
\end{lem}

\begin{lem}\label{l10} \mbox{\rm ( Critical case)}
\begin{enumerate}\item
Assume that $\a>1-1/p$, $1<p<\infty$  and $R>e$. There exist  positive numbers $C_6=C_6(\a,p,R )$ and  $L=L(\a,p,R )$ such  that  for  every $u\in C_c^\infty((0,1])$, we  have
\begin{align} \int_0^1\left ( |u'|^p -\Lambda_{\a,p} \left|\frac{u}{t}\right|^p \right ) t^{\a p}\,dt  + L|u(1)|^p 
\ge  C_6
\int_0^1 |u'|^p  \frac{t^{\a p}}{A_1^2}  \,dt. \label{3.9}
\end{align}
\item
Assume that $\a=1-1/p$, $1<p<\infty$  and  $R>e^e$. There exist  positive numbers $C_6=C_6(\a,p,R )$ and  $L=L(\a,p,R )$ such  that  for  every $u\in C_c^\infty((0,1])$, we  have
\begin{align} \int_0^1\left ( |u'|^p -\Lambda_{\a,p} \left|\frac{u}{t}\right|^p \frac{1}{A_1^p} \right ) t^{p-1}\,dt  + L|u(1)|^p
\ge  C_6
\int_0^1|u'|^p    \frac{t^{p-1}}{A_2^2}  \,dt.  \label{3.10}
\end{align}
\end{enumerate}

\end{lem}

\noindent{\bf Proof of Lemma \ref{l10}:} Admitting Lemma \ref{l7} for  the  moment, we  prove Lemma \ref{l10}.
Unfortunately  we can  not employ a counterpart of
Lemma \ref{l2}, hence  we  use  a  direct argument. We establish  (\ref{3.10}) (the  assertion 2 ) only because the  argument  for  
(\ref{3.9}) is quite  similar. 
We prepare  the following fundamental inequalities which are  established in \cite{anm}   as Lemma 2.1 for  $X>-1$. 
\begin{lem}\label{lem3.1} 
\begin{enumerate}\item
For $p\ge 2$ we have 
\begin{equation}
|1+X|^p-1-pX\ge  c(p)|X|^q, \quad \mbox{ for any } q\in [2,p]  \mbox{ and } X\in \mathbf {R}.
\label{3.0}\end{equation}
\item
For $1<p\le2$ and $M\ge1$, we  have 
\begin{equation} |1+X|^p-1-pX\ge  c(p)\begin{cases} &M^{p-2}X^2,\qquad |X|\le M,\\
& |X|^p,\qquad\qquad |X|\ge M.\end{cases}\label{3.1}\end{equation}
\end{enumerate}
Here $c(p)$ is a positive number independent of each $X$, $M\ge1$  and $q\in [2,p]$.
\end{lem}
{\noindent\bf Proof.} By  Taylar expansion  we have (\ref{3.0}) with $q=2$. 
For $p>1$, we  note that 
\begin{equation} \lim_{X\to 0} \frac{|1+X|^p-1-pX}{X^2}=\frac{p(p-1)}{2}, \quad   \lim_{|X|\to \infty} \frac{|1+X|^p-1-pX}{|X|^p} =1. \label{Add}
\end{equation}
Therefore  (\ref{3.0}) is valid for any $q\in [2,p]$ for a small $c(p)>0$.
If $X>-1$,  then (\ref{3.1}) also follows from Taylar expansion and  (\ref{Add}). 
If we choose $c(p)$ sufficiently small, then it
remains valid for $X\le -1$. 
\hfill$\Box$\par 

First  we  assume  that $p\ge 2$ and $ \a=1-1/p$. For  any $u\in C^1_c((0,1])$,
let  us  set $ u= A_1^{\a} h$, where $  A_1(t)= \log(R/t)$ and $ h\in C^1_c((0,1])$. 
Without a loss of  generality  we assume  $u\ge 0$. Letting $X=- {tA_1h' }/{(\a h)}\, (h\neq0); 0\, (h=0)$, we have 
\begin{align}|u'|^p t^{p-1}-\Lambda_{\a,p}u^p \frac{ t^{p-1}}{t^p A_1^p}
& =\Lambda_{\a,p} \frac{h^p}{tA_1}\left( |1+X|^p-1\right)   \label{3.14}\\
&\ge - (\Lambda_{\a,p})^{1-\frac 1p} (h^p)' +{c(p)} |h'|^p (A_1 t   )^{p-1}.\notag
\end{align}
 Here we  used (\ref{3.0}) with $q=p$.
On  the other hand we have 
\begin{align} |u'|^p  \frac{t^{p-1}}{A_2(t)^2} &= \a^p  \frac{h^p}{tA_1 A_2^2}|1+X|^p 
\le 2^p\a^p \frac{h^p}{tA_1 A_2^2} \left(1+|X|^p\right) \label{3.15} \\
&=2^p\a^p \frac{u^p}{tA_1^p A_2^2} + 2^p |h'|^p (tA_1)^{p-1} \frac{1}{A_2^2}.\notag
\end{align}
Here we used  a trivial  inequality:
$|1+X|^p \le  2^p (1+ |X|^p).
$
By integrating (\ref{3.14})  and (\ref{3.15}) on $(0,1)$ and employing Lemma  \ref{l7},  
 the desired inequality follows  for  a  sufficiently  small  constant $C_6>0$.\par
Secondly  we  assume  that $1<p<2$. 
If $|X|\ge  M$, then (\ref{3.15}) is  valid.
If $|X|\le M$, again from (\ref{3.15}) we immediately   have 
\begin{align} |u'|^p  \frac{t^{p-1}}{A_2(t)^2}  \le 2^p\a^p \left(1+M^p\right) \frac{h^p}{tA_1 A_2^2}= C(M) \frac{u^p}{tA_1^p A_2^2}  
 \qquad \left(  |X|\le M \right).\label{3.17} 
\end{align}
Thus we  have 
\begin{equation} |u'|^p  \frac{t^{p-1}}{A_2(t)^2} \le C(M)\frac{u^p}{tA_1^p A_2^2} + 2^p \chi_{|X|\ge M}(t)  |h'|^p (tA_1)^{p-1} \frac{1}{A_2^2}. \label{3.19}
\end{equation}
Here $\chi_S(t)$ is  a characteristic function of $S$. We have  (\ref{3.14})  provided  that $|X|\ge M$.
Since  $A_2^{-2} \le 1$, for  a  sufficiently  small $C_6$
the desired inequality (\ref{3.10})  follows  from (\ref{3.14}), (\ref{3.19}) and Lemma \ref{l7} (2).
\hfill $\Box$

\subsection{Proof of Theorem \ref{CT2} }
It follows  from (\ref{3.8}) and (\ref{3.10}) that we  have 
\begin{align} \int_0^1\left ( |u'|^p -\Lambda_{\a,p} \left|\frac{u}{t}\right|^p \frac{1}{A_1^p} \right )&  t^{p-1}\,dt  + L|u(1)|^p    \label{3.16} \\&
\ge  C_7
\int_0^1 \left(|u'|^p  +\Lambda_{\a,p} \left|\frac{u}{t}\right|^p \frac{ 1 }{A_1^p } \right)
 \frac{t^{p-1}}{A_2^2}  \,dt   \notag
\end{align}
for  every 
 $u\in C_c^\infty((0,1])$. Here $C_7= \displaystyle{{\min( C_5, \Lambda_{\a,p}C_6)}/{2}}$.
By  adding $$-C_0\int_0^1 \left|\frac{u}{t}\right|^p  \frac{1}{A_1(t)^p} \frac{1}{A_2(t)^2} t^{p-1}\,dt$$ to  the  both side
of  (\ref{3.10}) we  have  
\begin{align*}  \int_0^1  & \left( |u'|^p -\left|\frac{u}{t}\right|^p  \frac{1}{A_1(t)^p}\left(\Lambda_{\a,p}    + \frac{C_0 }{A_2(t)^2}\right) \right) t^{p-1}\,dt  +L|u(1)|^p \\ 
&\ge  C_7
 \int_0^1 \left( |u'|^p + \left|\frac{u}{t}\right|^p   \frac{1}{A_1(t)^p} \left( \Lambda_{\a,p}-\frac{C_0}{C_7} \right)\right) \frac{1}{A_1(t)^2} t^{p-1}\,dt,
\end{align*}
Now  we set $\, C_0={\Lambda_{\a,p} C_7}/{3} \, $ and  $\   C_1'= {C_7}/{3}.$
Assuming  that  $C_0\le \Lambda_{\a,p} (\log\log R)^2$,  we  have 
$ C_0/A_2(t)^2 \le  C_0 /( \log\log R)^2\le \Lambda_{\a,p}$.  Then
\begin{align*}  \int_0^1  & \left( |u'|^p -\left|\frac{u}{t}\right|^p  \frac{1}{A_1(t)^p}\left(\Lambda_{\a,p}    + \frac{C_0 }{A_2(t)^2}\right) \right) t^{p-1}\,dt  +L|u(1)|^p\\ 
&\ge  C_1'
 \int_0^1 \left( |u'|^p + \left|\frac{u}{t}\right|^p   \frac{1}{A_1(t)^p} \left( \Lambda_{\a,p}+ \frac{C_0 }{A_2(t)^2} \right)\right) \frac{1}{A_1(t)^2} t^{p-1}\,dt.
\end{align*}

By a  calculation  we see that   for  some $C(R)>0$,  $ t A_2^{2} \le C(R)\,  \mbox{ for  any } \,  t\in [0,1] . $
Therefore  the  desired  inequality holds  for  $C_1= C_1' C(R)^{-1}$.
 \hfill $\Box$

\section{ Proofs  of  Theorems \ref{NCT1}, \ref{CT3}  and  Corollaries \ref{NCC1}, \ref{CCT3}}
We first establish Theorem \ref{NCT1} using Theorem \ref{nct1}. 
 Theorem \ref{CT3} is   proved in  a quite similar  way using Theorem \ref{CT2}.
 Then we prove   Corollary  \ref{NCC1}  and   Corollary  \ref{CCT3}.
\par \medskip

\noindent{\bf  Proofs of Theorem \ref{NCT1} and  Theorem \ref{CT3}:} 
Let  us  prepare some notations and  fundamental facts.
Define $\Sigma = \partial \Omega$  and $\Sigma_t =\{ x\in  \Omega : \delta(x)= t\}$.
Since   $\Sigma$ is  is  of  class $C^2$, there exists an $\eta_0>0$  such  that we  have 
a $C^2$ diffeomorphism $ G: \Omega_\eta \mapsto (0,\eta)\times \Sigma$
for any  $\eta\in (0,\eta_0)$.  By $G^{-1}(t, \sigma) \, ((t,\sigma) \in (0,\eta_0) \times \Sigma)$ we denote the inverse of $G$. Let $H^{}_t$ denote the mapping 
$G^{-1}(t, \cdot) $ of $\Sigma$ onto $\Sigma_t$.
This mapping is also a $C^2$ diffeomorphism  and its Jacobian is close to $1$ in $(0,\eta_0)\times \Sigma$.
Therefore,
for every non-negative  continuous function $u$ on $\overline{\Omega_\eta}$ with $\eta\in (0,\eta_0)$  we  have 
\begin{align}& \int_{\Omega_\eta}u = \int _0^\eta \,dt \int_{\Sigma_t} u \,d\sigma_t =\int _0^\eta \,dt \int_\Sigma u(t, H_t(\sigma)) ( \rm{Jac}\, H_t )\,d\sigma,\\
& | \rm{Jac}\, H_t (\sigma)-1|\le ct, \quad \mbox{  for every   } (t,\sigma)\in (0,\eta_0)\times  \Sigma,
\end{align}
where $c$  is  a positive  constant independent of  each $(t,\sigma)$, $d\sigma$  and $d\sigma_t$  denote  surface elements on $\Sigma$  and $\Sigma_t$  respectively.
Further we  have 
 \begin{align} &  
  \int _0^\eta \,dt \int_\Sigma u(t, H_t(\sigma)) ( 1-ct)\,d\sigma \le  \int_{\Omega_\eta}u\le   \int _0^\eta \,dt \int_\Sigma u(t, H_t(\sigma)) ( 1+ct)\,d\sigma, \label{4.1}\\ 
 & 
\int_\Sigma u(\eta, H_\eta(\sigma)) ( 1-c\eta)\,d\sigma \le 
\int_{\Sigma_\eta}u\,d\sigma_\eta\le   \int_\Sigma u(\eta, H_\eta(\sigma)) ( 1+c\eta)\,d\sigma. \label{4.2}
\end{align}

Then we immediately  have 
\begin{align}
  \int_{\Sigma} \,d\sigma \int_0^\eta \left |\frac{\partial u}{\partial t} \right |^p (1-ct) t^{\a p} \,dt &\le \int_{\Omega_\eta}|\nabla u|^p \delta^{\a p}
  \notag\\
 \int_{\Sigma} \,d\sigma \int_0^\eta \left |\frac{ u}{t} \right |^p (1-ct) t^{\a p} \,dt \le  \int_{\Omega_\eta}| u|^p \delta^{p(\a -1)}&\le  \int_{\Sigma} \,d\sigma \int_0^\eta \left |\frac{ u}{t} \right |^p (1+ct) t^{\a p} \,dt. \notag
 \end{align}
 \noindent{\bf  Proof  of (\ref{2.6}):} Under these consideration, 
 (\ref{2.6}) is  reduced  to  one-dimensional Hardy's inequality. Setting $v(t) =u(t,\sigma)$  and $v'= \partial u/\partial t$ we  have
 
 \begin{align}  \int_0^\eta  & \left( |v'|^p -\left|\frac{v}{t}\right|^p \left(\Lambda_{\a,p}    + \frac{C_2 }{A_1(t)^2}\right)\right) t^{\a p}\,dt  \notag \\ &\ge  c
 \int_0^\eta \left( |v'|^p + \left|\frac{v}{t}\right|^p \left(\Lambda_{\a,p} + \frac{C_2 }{A_1(t)^2}\right) \right) t^{\a p+1}\,dt+ L |v(\eta)|^p \eta^{\a p}(1+c\eta).
\end{align}
By a change of  variable $t=s \eta$, putting $w(s)= v( s\eta)$ with $v\in C^1_c((0,1])$,

 \begin{align}&  \int_0^1   \left( |w'|^p -\left|\frac{w}{s}\right|^p \left(\Lambda_{\a,p}    + \frac{C_2 }{A_1(s\eta)^2}\right)\right)s^{\a p}\,ds  \\ &\ge  c\eta
 \int_0^1 \left( |w'|^p + \left|\frac{w}{s}\right|^p \left(\Lambda_{\a,p} + \frac{C_2 }{A_1(s\eta)^2}\right) \right) s^{\a p+1}\,ds+ L|w(1)|^p\eta^{p-1}(1+c\eta).\notag
\end{align}

 On  the  other hand,
by Theorem \ref{nct1}  with  $R$ changed to   $R/\eta$, we have, for every  $w\in  C^1_c((0,1])$,

\begin{align}  \int_0^1  & \left( |w'|^p -\left|\frac{w}{t}\right|^p \left(\Lambda_{\a,p}    + \frac{C_0 }{A_1(t\eta)^2}\right)\right) t^{\a p}\,dt  \\ &\ge  C_1
 \int_0^1 \left( |w'|^p + \left|\frac{w}{t}\right|^p \left(\Lambda_{\a,p} + \frac{C_0 }{A_1(t\eta)^2}\right) \right) t^{\a p+1}\,dt + L |w(1)|^p,\notag 
\end{align}
where $C_0$  and  $C_1$  may  depend  on $\eta$  but   independent  of  each  function $v$.
Now  we take  $\eta$ and  $C_2$ so  that $C_1/\eta> c$,  $C_2 \le C_0$ and  $\eta^{p-1}(1+c\eta) <1$ respectively.  Since $w$ is an arbitrary function in $ C^1_c((0,\eta])$,   we get  (\ref{2.6}). 
\hfill$\qed$
\par\bigskip
 \noindent{\bf  Proof  of (\ref{2.11}) and (\ref{2.12}):}
 In  parallel to  the  verification of (\ref{2.6}),  (\ref{2.11}) and (\ref{2.12}) can be proved using (\ref{c1}) and  (\ref{2.10}) together  with (\ref{4.2}).
 Hence  we  omit the detail.\hfill$\qed$
\par\bigskip
 
\noindent{\bf  Proof of Corollary \ref{NCC1}:} 
Assume on  the contrary
 that Hardy inequality (\ref{2.7}) does  not  hold. Then there exists a sequence of  functions $\{u_k\} \subset W^{1,p}_{\alpha, 0}(\Omega)$ such  that
\begin{equation}  \lim_{k\to \infty}\int_{\Omega} |\nabla u_k|^p\delta^{\alpha p}\,dx =0,\quad
 \int_{\Omega} {|u_k|^p}\delta^{p(\alpha-1)}\, dx =1 \quad (k=1,2,\cdots).
\end{equation}

By Theorem {\ref{NCT1}} we  have 
\begin{align*} 
&\int_{\Omega} |\nabla u_k|^p\delta^{\alpha p}\,dx =  \int_{\Omega_\eta} |\nabla u_k|^p\delta^{\alpha p}\,dx + \int_{\Omega\setminus \Omega_\eta} |\nabla u_k|^p\delta^{\alpha p}\,dx\\
&\ge  \Lambda_{\alpha ,p}\left( 1-  \int_{\Omega \setminus \Omega_\eta} {|u_k|^p}\delta^{p(\alpha-1)}\, dx\right)
+L \int_{\Sigma_\eta} |u_k|^p\delta^{\a p} 
+ \int_{\Omega\setminus \Omega_\eta} |\nabla u_k|^p\delta^{\alpha p}\,dx.
\end{align*}
Since  $\delta \ge \eta$ in $\Omega\setminus \Omega_\eta$,  by  the  standard argument we have $u_k\to C \,( constant)$ in $W^{1,p}(\Omega\setminus \Omega_\eta)$ as $k\to \infty$. Since  $L>0$, we  have $C=0$.
Hence  
 we see
$ 0 \ge \Lambda_{\alpha ,p}$, and   we  reach  to   a contradiction.  \hfill $\Box$
\par
\noindent{\bf  Proof of Corollary \ref{CCT3}:} 
Assume on  the contrary
that Hardy inequality (\ref{2.14}) does  not  hold. Then there exists a sequence of  functions $\{u_k\} \subset W^{1,p}_{\alpha, 0}(\Omega)$ such  that
\begin{equation} \begin{cases}& \lim_{k\to \infty}\int_{\Omega} |\nabla u_k|^p\delta^{\alpha p}\,dx =0, \,
\lim_{k\to\infty}\int_{\Sigma_\eta} |u|^p\delta^{\a p} =0,\\
& \int_{\Omega} {|u_k|^p}\delta^{p(\alpha-1)}\, dx =1 \quad (k=1,2,\cdots).
\end{cases}
\end{equation}
By Theorem {\ref{CT3}} we  have 
\begin{align*} 
&\int_{\Omega} |\nabla u_k|^p\delta^{\alpha p}\,dx +L \int_{\Sigma_\eta} |u_k|^p\delta^{\a p} \\
&\ge  \Lambda_{\alpha ,p}\left( 1-  \int_{\Omega \setminus \Omega_\eta} {|u_k|^p}\delta^{p(\alpha-1)}\, dx\right)
+ \int_{\Omega\setminus \Omega_\eta} |\nabla u_k|^p\delta^{\alpha p}\,dx.
\end{align*}
Since  $\delta \ge \eta$ in $\Omega\setminus \Omega_\eta$,   as  before we have $u_k\to 0$ in $W^{1,p}(\Omega\setminus \Omega_\eta)$ as $k\to \infty$. Hence  
 we have  
$ 0 \ge \Lambda_{\alpha ,p}$, and   we  get   a contradiction. \hfill $\Box$

\section{Proof of Theorem \ref{NCT2} }
It suffices  to  show  the  implication $1\to 2$. Since 
 $A_1(\delta)^{-1}  \le 1$  in $\Omega$  and  the  trace operator $T: W^{1,p}_{0,\a}(\Omega_\eta^c) \mapsto 
 L^p(\Sigma_\eta; \delta^{\a p} )$ is  continuous for  a  small  $\eta>0$, one  can assume  that   $C_2= L=0$.
Now, we assume on  the contrary   that
there exists a sequence of  functions $\{u_k\} \subset W^{1,p}_{\alpha, 0}(\Omega)$ such  that
\begin{equation} \label{3.17bis} \lim_{k\to \infty}\int_{\Omega_\eta} |\nabla u_k|^p\delta^{\alpha p}\,dx =0,\quad
 \int_{\Omega_\eta} {|u_k|^p}\delta^{p(\alpha-1)}\, dx =1 \quad (k=1,2,\cdots).
\end{equation}
Here we  prepare a lemma on   extension:

\begin{lem}\label{l11}\mbox{\rm ( Extension )} For any $\eta >0 $ there exists an extension operator  $E=E(\eta) : W^{1,p}_{\a,0}(\Omega_\eta)\mapsto  W^{1,p}_{\a,0}(\Omega)$  such  that:
\begin{enumerate}\item $E (u)= u \quad $ a.e. in $\Omega_\eta$
\item There exists some  positive  number $C=C(\eta)$  such that  for  any $u\in W^{1,p}_{\a,0} (\Omega_\eta)$,
$$ || |\nabla E( u)|||_{L^p(\Omega, \delta^{\a p})} \le  C\Big(  |||\nabla u|||_{L^p (\Omega_{\eta/2},\delta^{\a p} )}+ ||u||_{W^{1,p}_{\a,0} (\Omega_\eta \setminus \Omega_{\eta/2})}\Big).$$
\end{enumerate}

\end{lem}
Admitting this  for  the  moment, we prove  Theorem \ref{NCT2}.
Let $v_k= E(u_k) \in W^{1,p}_{\a,0}(\Omega)$ for $k=1,2,\ldots$. 
It  follows  from (\ref{3.17}), the  assumption 1 and the property  of $E$ that 
$v_k$ becomes  a Cauchy sequence  and  $v_k\to v $ in $ W^{1,p}_{\a,0}(\Omega)$ for  some $v\in  W^{1,p}_{\a,0}(\Omega)$ as $k\to\infty$.
On  the  other hand 
 by choosing  a subsequence if  necessary, we  see  that 
$u_k\to c$ a.e. in $\Omega_\eta$   for  some  constant  $c$ as $k\to\infty$. Then, by the  assumption 1
 \begin{align*} 1 &\le  \int_\Omega |E(u_k)|^p \delta^{p(\a-1)}\le  \gamma^{-1}
 \int_\Omega |\nabla E(u_k)|^p\delta^{\a p} \\ 
& \le    \gamma^{-1}C\Big(  |||\nabla u_k|^p||_{L^p (\Omega_{\eta/2},\delta^{\a p} )}+ ||u_k||_{W^{1,p}_{\a,0} (\Omega_\eta \setminus \Omega_{\eta/2})}\Big)<\infty.
 \end{align*}
Since $(\a-1)p<-1$,  we  have $c=0$.
Thus  $u_k \to 0 $  in $L^p(\Omega_\eta\setminus \Omega_{\eta/2})$.
Thus  we see that $ ||u_k||_{W^{1,p}_{\a,0} (\Omega_\eta \setminus \Omega_{\eta/2})}\to  0 $  as  $k\to\infty$. From  this  together  with (\ref{2.6}) we  have  a  contradiction.
 \hfill $\Box$
\par\medskip
 \noindent{\bf Proof  of  Lemma \ref{l11}: } Since  $\delta$ is  Lipschitz continuous, 
we  see  that $\partial \Omega_\eta$  and $\partial \Omega_{\eta/2}$ are Lipschitz   compact manifolds. By  the  standard theory  we  have  an extension operator $\tilde{E} : W^{1,p} (\Omega_\eta \setminus \Omega_{\eta/2})\mapsto
 W^{1,p}(\Omega\setminus \Omega_{\eta/2})$ such  that $\tilde E(u)= u$  a.e. in $\Omega_\eta \setminus \Omega_{\eta/2}$, and 
$$ || |\nabla \tilde{E}( u)|^p||_{L^p(\Omega\setminus \Omega_{\eta/2}, \delta^{\a p})} \le  C(\eta)  ||u||_{W^{1,p}_{\a,0} (\Omega_\eta \setminus \Omega_{\eta/2})}.$$
 Define  for  $u\in W^{1,p}_{\alpha, 0}(\Omega_\eta)$
  \begin{equation}  E(u) = u \,\, (x\in \Omega_{\eta/2},);  \quad  \tilde{E}(u) \,\,  (x  \in  \Omega\setminus \Omega_{\eta/2} ).
 \end{equation}
 Then the  assertion follows.\hfill$\Box$
 
 \section{Proofs  of  Propositions \ref{ct1}  and  \ref{CT1}}
Proposition  \ref{ct1} is    known in a more general  fashion. In  fact  
a variant  is seen in Maz'ya   \cite{Ma} ( Lemma 2, p144). For  the sake  of  reader's convenience we  give a direct verification. 
We note  that 
Proposition \ref{CT1} is  a   consequence  of  Proposition  \ref{ct1}.
\par
\noindent{\bf Proof of Proposition \ref{ct1}:} 
First  we assume  that  $ \a> 1-\frac 1p$. Define  for $\varepsilon \in (0,1)$,  
$\,  u_\varepsilon=  t/{\varepsilon}\, (0\le t\le \varepsilon);\,
 1\,  ( t\ge \varepsilon).$ 
Then we immediately  have $u_\varepsilon(0)=0, u_\varepsilon(1)=1$  and 
$ \int _0^1|u _\varepsilon'|^p t^{\a p} \,dt
\to 0\,\mbox{ as } \varepsilon\downarrow 0.$
By using $C^1$ approximation of each $u_\varepsilon$, the  assertion is  proved.
Further  we note  that
$ \int_0^1 |u_\varepsilon|^p t^{(\a-1)p}\,dt 
\to 1/ ({\a p-p+1})>0
\,  \mbox{ as } \varepsilon\downarrow 0.$
In  the  critical  case, define for $\varepsilon \in (0,1/2)$
 \begin{equation} u_\varepsilon= \,
0 \,( 0\le t\le \varepsilon) ;  \quad \frac {A_1(\varepsilon)-A_1(t)  }{A_1(\varepsilon)-A_1(1/2) }\, \,(\varepsilon \le t\le 1/2);\quad  1\, \,(1/2 \le t \le 1). \end{equation}
Then 
$ \int _0^1|u_\varepsilon '|^p t^{p-1} \,dt=( A_1(\varepsilon)-A_1(1/2))^{1-p} \to 0 \, \mbox{ as } \varepsilon\downarrow 0.$ 
On  the other hand  we  have  $u_\varepsilon(0)=0, u_\varepsilon(1)=1$ and  hence the  assertion is now  clear. 
Further we  note  that 
$$\int_\varepsilon^1 |u_\varepsilon|^p \frac{1}{t A_1(t)^p}\,dt\ge  \frac{A_1(1)^{1-p}- A_1(1/2)^{1-p}}{p-1}  >0
\quad  \mbox{ as } \varepsilon\downarrow 0.$$
\qed

\noindent{\bf  Proof  of  Proposition \ref{CT1}:}
We give  a  proof  when $\a>1-\frac 1p$, because  the  argument is quite  similar in  the rest  of  the  case.
If  a positive  number $\eta_0$  is sufficiently  small, then one  can  assume  that  $\delta\in C^2(\Omega_{\eta_0})$  and 
a manifolds $\{ x\in \Omega;  \delta =\eta\}$ is of  $C^2$  class for  $\eta \in (0,\eta_0]$. 
By  virtue  of  (\ref{4.1}) we  have
$$ \int_{\Omega_\eta}| u|^p \delta^{p(\a -1)}\le  \int_{\Sigma} \,d\sigma \int_0^\eta \left |\frac{ u}{t} \right |^p (1+ct) t^{\a p} \,dt,$$
hence  the  assertion follows  from  Proposition \ref{ct1}.\qed

\section{Appendix; Proofs of Lemma \ref{l3} and Lemma \ref{l7} }

\subsection{Preliminary  }
In this section we prepare  a series of one dimensional weighted Hardy's inequalities.
The followings are  given in  \cite{ah0} as Lemma {3.1} and Lemma {3.4} respectively.
\begin{lem}\label{lem3.3}
 Assume that  $R>e$. Then, for any $h\in  C^1_c((0,1])$
 we have
\begin{equation}
\int_0^1|h^\prime(t)|^2 tdt\geq\frac14
\int_0^1|h(t)|^2A_1(t)^{-2}\frac{dt}{t} -\frac12 A_1(1)^{-1}h(1)^2.\label{6.3}
\end{equation}
\end{lem}
{\noindent\bf  Proof.}  Let $h(t)=A_1(t)^\frac12 w(t)$. Then   we have
\begin{align}
|h^\prime(t)|^2t=\frac{t}{A_1(t)}\left(-\frac{1}{2t}w(t)+w^\prime(r)A_1(t)\right)^2
&\geq \frac{|h(t)|^2}{4t A_1(t)^{2}}-
\frac{1}{2}\left(\frac{d}{dt}w^2(t)\right).
\label{6.5}\end{align}
Since $w(0)=0$, 
we have (\ref{6.3}) and 
the rest of  the proof is  clear.\hfill$\Box$\par
\begin{lem}\label{lem3.4}
 Assume that $R>e^e$. Then,   for any $h\in C^1_c((0,1])$  we have 
\begin{equation}
\int_0^1|h^\prime(t)|^2 t  A_1(t) \,dr\geq \frac14
\int_0^1
\frac{|h(t)|^2}{  A_1(t)\cdot A_2(t)^{2}}\frac{dt}{t} -\frac12 A_2(1)^{-1}h(1)^2.\label{6.6}
\end{equation}
\end{lem}
{\noindent\bf  Proof.}  
For $h(t)= A_2(t)^\frac12 w(t)$, we  have in a similar way 
\begin{align}
|h^\prime(t)|^2 t A_1(t)
\geq\frac{1}{4}\frac{|h(t)|^2}{ A_1(t)
 A_2(t)^{2}}-
\frac{1}{2}\left(\frac{d}{dt}w^2(t)\right).
\end{align}
Then  the  rest of the proof is  clear.\hfill$\Box$\par
\begin{df}\label{df3.1} A function $\varphi\in C^1([0,1])$ is said to belong to $G([0,1])$ 
if and only if 
\begin{equation} \varphi(0)= 0, \quad   \varphi'(0)\neq 0 \quad  \text{ and  }  \quad \varphi'(1)=0.  \label{case1}\end{equation}
\end{df}  
\begin{df}\label{df3.2} For $\varphi \in G([0,1])$ and $M>1$ we define  three subsets of $[0,1]$ as
follows:
\begin{equation}\begin{cases}&A(\varphi,M)=\left\{ t\in [0,1] \,| \,|\varphi'(t)|\le
M\frac{|\varphi(t)|}{t}
\right\},\\
& B(\varphi,M)=\left\{ t\in [0,1] \,| \,|\varphi'(t)|> M\frac{|\varphi(t)|}{t} \right\},\\
& C(\varphi,M)=\left\{ t\in [0,1] \,| \,|\varphi'(t)|= M\frac{|\varphi(t)|}{t} \right\}.\end{cases}
\label{6.10}\end{equation}
\end{df}
Clearly  $[0,1] =A(\varphi,M)\cup B(\varphi,M)$. From (\ref{case1})  we  see  $0, 1\in A(\varphi,M)$. 
We    note  that
the set $ C(\varphi,M)$ coincides with the set of
critical points of $ \log (|\varphi|r^{\pm M})$. 
By  a  standard argument we  have  the  following  approximation lemma ( cf.   Lemma 3.5 in \cite{ah0}).

\begin{lem}\label{lem3.5}Let $M>1$  and $\varphi\in G([0,1])\cap C^2([0,1])$. Assume that $\varphi\ge 0$. Then there
exists a sequence  of functions
$\varphi_k\in G([0,1])\cap C^2([0,1])$ such that $\varphi_k>0$ in $(0,1)$, $\varphi_k\to \varphi$ in
$C^1([0,1])$ as
$k\to+\infty$ and 
$C(\varphi_k,M) $ consists of finite points for any
$k$.
\end{lem}
We  prepare some estimates for the proofs of  Lemma \ref{l3} and Lemma \ref{l7}. 
\begin{lem}\label{lem3.6} Assume that $1<p<2$ and $R>e$. Then for any $\varepsilon>0$ there is a  positive
number $M$ such that we have for any $\varphi\in G([0,1])$
\begin{equation}\int_{B(\varphi,M)}\frac{|\varphi|| \varphi'|}{  A_1(t)}\,dt\le \varepsilon 
\int_{B(\varphi,M)}|\varphi|^{2-p}|\varphi'|^p t^{p-1}\,dt. \label{6.14}
\end{equation}
\end{lem}
{\noindent\bf Proof.} 
We may assume that $\varphi>0$. Then by the definition we have $t{|\varphi'|}/{\varphi}>M
$ on $B(\varphi,M)$. Hence we immediately have 
\begin{equation} \varphi^{2-p}|\varphi'|^p t^{p-1}= \varphi| \varphi'| \cdot
\left(t\frac{|\varphi'|}{\varphi}\right)^{p-1}\ge  M^{p-1}\varphi| \varphi'|, \qquad\text{ on }
B(\varphi,M).\end{equation} 
Therefore it suffices to choose $M$ so that
$ M^{1-p}(\log R)^{-1}\le \varepsilon$.
\hfill$\Box$\par
\begin{lem} \label{lem3.7} Assume that $1<p<2$ and $R>e$. Then we have  for any $\varphi\in G([0,1])$
\begin{align}\int_{A(\varphi,M)}|\varphi'(t)|^2t\,dt &\ge
\frac14\int_{A(\varphi,M)}\frac{|\varphi|^{2}}{t A_1(t)^2}\,dt   -\frac 12 \frac{\varphi(1)^2}{A_1(1)}\label{6.17}\\&+\frac12 
\int_{B(\varphi,M)}\frac{|\varphi|^{2}}{t A_1(t)^2}\,dt-
\int_{B(\varphi,M)}\frac{|\varphi||\varphi'|}{ A_1(t)}\,dt.\notag
\end{align}
\end{lem}
{\noindent\bf Proof.} By Lemma \ref{lem3.5}
we can assume that  $C(\varphi,M)$ consists of
finitely many points.  Recall that $0,1\in A(\varphi,M)$. 
From  the  argument of Lemma \ref{lem3.3} we  have
\begin{align}\int_{A(\varphi,M)}|&\varphi'(t)|^2t\,dt\ge
\frac14\int_{A(\varphi,M)}\frac{|\varphi|^{2}}{t A_1(t)^2}\,dt
-\frac12 \int_{A(\varphi, M)}\frac{d}{dt}\left(\frac{\varphi(t)^2}{  A_1(t)}\right)\,dt
\label{6.19}\\
&=\frac14\int_{A(\varphi,M)}\frac{|\varphi|^{2}}{t A_1(t)^2}\,dt
+\frac12 \int_{B(\varphi,M)}\frac{d}{dt}\left(\frac{\varphi(t)^2}{  A_1(t)}\right)\,dt -\frac 12 \frac{\varphi(1)^2}{A_1(1)}. \notag
\end{align}
Thus we have  the desired estimate.
\hfill$\Box$\par
In a quite  similar way we have 
\begin{lem}\label{lem3.8} Assume that $1<p<2$ and $ R>e^e $. Then we have  for any $\varphi\in G([0,1])$
\begin{align}\int_{A(\varphi,M)}|\varphi'(t)|^2r   A_1(t)\,dt &\ge
\frac14\int_{A(\varphi,M)}\frac{|\varphi|^{2}}
{t A_1(t) A_2(t)^2}\,dt     -\frac 12 \frac{\varphi(1)^2}{A_2(1)}\label{6.20}\\&+\frac12 
\int_{B(\varphi,M)}\frac{|\varphi|^{2}}{t  A_1(t) A_2(t)^2}\,dt-
\int_{B(\varphi,M)}\frac{|\varphi||\varphi'|}{A_2(t)}\,dt. \notag
\end{align}

\end{lem}

\subsection{Proof  Lemma \ref{l3}  }
Assume that $ \a<1-1/p$.  
For $u\in C^1_c((0,1])$, we  define
\begin{equation}
u(t)=h(t)t^{\beta},\quad \beta= 1-\frac{1}{p} -\a,\quad( \beta^p= \Lambda_{\a,p}).\label{6.21}
\end{equation}
Without  the loss of generarity  we assume that $u\ge0$ in $(0,1)$, then we have 
\begin{equation}\begin{split}
\int_{0}^1|u'|^p t^{\a p}\, dt-&\Lambda_{\a,p}\int_{0}^1\frac{|u|^p}{t^{p}}t^{\a p}\,dt 
=\Lambda_{\a,p}\int_0^1 h^p
\left\{\left|1+\frac{rh'}{\beta h}\right|^p-1\right\}\frac{dt}{t}.\end{split}\label{6.23}
\end{equation}
For the moment we assume $p\ge 2$.
By the fundamental inequality (\ref{3.0}) with $q=2$, we obtain
\begin{equation}\begin{split} &\text{(R.H.S.) of (\ref{6.23})}
\ge \frac{\Lambda_{\a,p}}{\beta}\int_0^1p h^{p-1}h'\,dt+
c(p)\frac{\Lambda_{\a,p}}{\beta^2} \int_0^1h^{p-2}(h')^2t\,dt\\
&=\beta^{p-1} h(1)^p+
c(p) \beta^{p-2}\frac{4}{p^2}\int_0^1\left|\left(h^{\frac{p}{2}}(t)\right)^\prime
\right|^2tdt.\, \quad(\text{ Note  that  }
h(0)=0.)
\end{split} \label{spit}\end{equation}
Using Lemma \ref{lem3.3} we get
\begin{align*}
\int_0^1\left|\left(h^{\frac{p}{2}}(t)\right)^\prime\right|^2tdt 
\ge \frac{1}{4}\int_{0}^1\frac{|u(t)|^p}{t^{p}} A_1(r)^{-2}t^{\a p}dt 
 -\frac 12  A_1(1)^{-1} h(1)^p.
\end{align*}
Combining this with (\ref{spit}), we get the  inequality (\ref{3.3}) with $C_3=c(p) \beta^{p-2}/{p^2}$ and $L=\beta^{p-1}-2C_3 A_1(1)^{-1}
$, making $c(p)$ smaller if  necessary.
\par We proceed  to  the  case that $1<p<2$. 
For
 $u\in C^1_c((0,1])$, we retain the  notation (\ref{6.21}). Suppose that $M$ is sufficiently large.
In   Definition \ref{df3.2} 
we replace
$\varphi$ and $M$ by $h$ and $\beta M$ respectively, and assume that $h\in G([0,1])$ again. 
Lemma \ref{lem3.1} (2) implies
\begin{align}
&\int_{0}^1|u'|^pt^{\a p}\, dt-\Lambda_{\a,p}\int_{0}^1\frac{|v|^p}{t^{p}}t^{\a p}\,dt
=\Lambda_{\a,p}\int_0^1 h^p(t)
\left\{\left|1+\frac{th'}{\beta h}\right|^p-1\right\}\frac{dt}{t}\label{6.26} 
\\ 
&\ge  \beta^{p-1}h(1)^p+
 \frac{4  c(p) (M\beta)^{p-2}}{p^2}\int_{ A(h,\beta M)}((h^{\frac{p}{2}})')^2t\,dt
+c(p) \int_{ B(h,\beta M)}|h'|^{p} 
t^{p-1}\,dt.\notag
\end{align}
Using Lemma \ref{lem3.7} with  $A(h,\beta M)=A(h^{\frac{p}{2}},{p\beta M}/{2})$  and 
$B(h,\beta M)=B(h^{\frac{p}{2}},{p\beta M}/{2})$,  
\begin{align}
&\int_{ A(h,\beta M)}((h^{\frac{p}{2}})')^2t\,dt
\label{6.27}\\&\ge
\frac14\int_{ A(h,\beta M)}\frac{h^p}{t A_1(t)^2}\,dt -\frac12 \frac{h(1)^p}{ A_1(1)}+\frac12
\int_{ B(h,\beta M)}\frac{h^p}{t A_1(t)^2}\,dt
-\frac{p}{2}\int_{ B(h,\beta M)} \frac{h^{p-1}|h'|}{ A_1(t)}\,dt.\notag
\end{align}
We can estimate the last term to obtain 
\begin{equation}\frac{p}{2}\int_{ B(h,\beta M)} \frac{h^{p-1}|h'|}{ A_1(t)}\,dt\le 
\frac{p}{2}\frac{1}{(\beta M)^{p-1}\log R}
\int_{ B(h,\beta M)} |h'|^p t^{p-1}\,dt.\label{6.28}
\end{equation}
Here we simply used the fact that   $t|h'|>\beta M h$ holds on the set $B(h,\beta M)$.
Combining this with (\ref{6.26}) and (\ref{6.27}) for sufficiently  large $M$, we have the desired inequality.\hfill $\Box$

\subsection{Proof of Lemma \ref{l7} }
We  treat 
 the  case $\a=1-{1}/{p}$ only, because the    argument for  $\a>1-1/p$ is  similar to the  previous subsection. 
 For $u\in C^1_c((0,1])$ we  define \begin{equation}u(t)=  A_1(t)^\beta h(t), \quad \beta=
1-\frac{1}{p}.\label{6.29}
\end{equation}  Without  loss of  generality  we assume  $ u,h\ge  0$. 
 First   we assume $p\ge 2$.
Using (\ref{3.0}) with $q=2$ and 
$X=- \beta^{-1} t A_1(t) h'(t)  h(t)^{-1}$,
we obtain
\begin{equation}
\begin{split}
&\int_{0}^1| u'|^pt^{p-1}\, dt-\beta^p\int_{0}^1\frac{|u(t)|^p }{t A_1(t)^{p}}\,dt \\&
=\beta^p\int_0^1  \frac{h(t)^p}{tA_1(t)} \left( \left|1- \frac{ t A_1(t) h'(t)}{\beta
h(t)}\right|^p-1
\right)\,dt\\&
\ge 
-\beta^{p-1}  h(1)^p + \frac{4 c(p){\beta^{p-2}}}{p^2 }\int_0^1\left|\left(h(t)^{\frac{p}{2}}(t)\right)^\prime
\right|^2t A_1(t)\,dt.\\
\end{split}\label{6.31}
\end{equation}
Using Lemma \ref{lem3.4}  we get
\begin{equation}
\int_0^1\left|\left(h(t)^{\frac{p}{2}}\right)^\prime\right|^2t 
A_1(t)\,dt   
\ge \frac{1}{4}\int_{0}^1\frac{|u(t)|^p}{t A_1(t)^pA_2(t)^2}dt  -\frac 12 A_2^{-1}(1)h(1)^p.
\label{6.32}\end{equation}
Combining this with (\ref{6.31}) we get the desired inequality where $C_5=c(p)\beta^{p-2}p^{-2}$ and 
$L= A_1^{1-p}(\beta^{p-1}+2C_5 A_2^{-1})$.
Then we proceed  to    the  case that $1<p<2$. 
Since   the  argument is
quite similar, 
  we  give  a sketch of  proof.
Suppose that $M$ is sufficiently large. We retain the notation (\ref{6.29}) and  modify Definition \ref{df3.2} as  follows:

\begin{df}\label{df4.1}For $\varphi \in G([0,1])$ and $M>1$ we define  three subsets of $[0,1]$ as
follows:
\begin{equation}\begin{cases}& A(\varphi,M)=\left\{ t\in [0,1] \,| \,|\varphi'(t)|\le
M\frac{|\varphi(t)|}{tA_1(t)}
\right\},\\
& B(\varphi,M)=\left\{ t\in [0,1] \,| \,|\varphi'(t)|> M\frac{|\varphi(t)|}{tA_1(t)} \right\},
\\
& C(\varphi,M)=\left\{ t\in [0,1] \,| \,|\varphi'(t)|= M\frac{|\varphi(t)|}{tA_1(t)} \right\}.
\end{cases}\label{6.33}\end{equation}
\end{df}
Again we  replace $\varphi$ and $M$ by $h$ and $\beta M$ respectively  and  assume  
 $h\ge 0 $ in $ (0,1)  $.
By Lemma {\ref{3.8}}, the  first line of (\ref{6.31}) is  estimated   from  below by  the  following:
\begin{align}
&
- {\beta^{p-1}}\int_0^1p h^{p-1}h'\,dt
+
\beta^p\int_{ A(h,\beta M)}h^{p} c(p)  M^{p-2}
\left(\frac{tA_1(t)h'}{\beta h}\right)^2\frac {A_1(t)^{-1}} t\,dt \label{6.34}
\\&
+\beta^p \int_{ B(h,\beta M)}h^{p} c(p) 
\left|\frac{tA_1(t)h'}{\beta h}\right|^p\frac{A_1(t)^{-1}}t\,dt
\notag
\\&
= -\beta^{p-1} h(1)^p+c(p)  (M\beta)^{p-2}\int_{ A(h,\beta M)}h^{p-2}|h'|^2tA_1(t)\,dt 
\notag
\\&
+c(p) \int_{ B(h,\beta M)}|h'|^{p} 
t^{p-1}A_1(t)^{p-1}\,dt.
\notag
\end{align}
Here we  note  that  $A(h,\beta M)=A(h^{\frac{p}{2}},{p}\beta M/{2})$ and $B(h,\beta M)=B(h^{\frac{p}{2}},{p}\beta
M/{2})$. Then applying Lemma {\ref{lem3.8}} for  $\varphi =h^{\frac{p}{2}}$,  $A(h^{\frac{p}{2}},{p}\beta M/{2})$ and  $B(h^{\frac{p}{2}},{p}\beta
M/{2})$ we  have 
\begin{align}
\int_{ A(h,\beta M)}&h^{p-2}(h')^2tA_1(t)\,dt= \frac{4}{p^2}\int_{ A(h,\beta
M)}((h^{\frac{p}{2}})')^2tA_1(t)\,dt\label{6.35}\\&\ge
\frac{4}{p^2}\bigg(
\frac14\int_{A(h,\beta M)}\frac{h(t)^p}
{tA_1(t) A_2(t)^2}\,dt -\frac12  A_2(1)^{-1} h(1)^p \notag\\&+\frac12 
\int_{B(h,\beta M)}\frac{h(t)^p}{tA_1(t) A_2(t)^2}\,dt-\frac{p}{2}
\int_{B(h,\beta M)}\frac{h(t)^{p-1}|h'(t)|}{A_2(t)}\,dt\notag\bigg).\notag
\end{align}
From an easy variant of Lemma \ref{lem3.6} we can estimate the last term to obtain 
\begin{equation}\frac{p}{2}\int_{ B(h,\beta M)} \frac{h^{p-1}|h'|}{A_2(t)}\,dt\le 
\frac{p}{2}\frac{1}{(\beta M)^{p-1}\log(\log R)}
\int_{ B(h,\beta M)} |h'|^p A_1(t)^{p-1}t^{p-1}\,dt.\label{6.36}
\end{equation}
Here we simply used the fact that   $tA_1(t)|h'|>\beta M h$ holds on the set $B(h,\beta M)$.
Combining this with (\ref{6.34}) and (\ref{6.35}), for a  large $M$, we have the desired inequality. \hfill $\Box$

\bigskip
Xiaojing Liu:  Department Math., Ibaraki University, Mito Ibaraki 310-8512, Japan; 
eleven11qq@163.com\\ \\

Hiroshi Ando: 
 Department Math., Ibaraki University, Mito Ibaraki 310-8512, Japan;  hiroshi.ando.math@vc.ibaraki.ac.jp  \\ \\
 Toshio Horiuchi: Department Math., Ibaraki University, Mito Ibaraki 310-8512, Japan; toshio.horiuchi.math@vc.ibaraki.ac.jp 


\begin{thebibliography}{BVxx}

\bibitem[1]{anm} Adimurthi, N. Nirmalendu, M. Chaudhuri, and Mythily Ramaswamy, An improved Hardy-Sobolev
inequality and its application, {\it Proceedings of the American Mathematical Society}, {\bf Vol. 130},
No. 2, 2001, pp. 489-505.

\bibitem[2]{ah0}  H. Ando, T. Horiuchi, Missing terms in the weighted  Hardy-Sobolev inequalities and its  application, 
 {\it Kyoto Journal of Mathematics,} {\bf Vol. 52}, No. 4,  (2012), pp. 759-796. 

\bibitem[3]{ah2}  H. Ando, T. Horiuchi,  Weighted Hardy's inequalities and the variational problem 
with compact perturbations, {\it Mathematical Journal of Ibaraki University,}  {\bf Vol. 52},  (2020), pp. 15-26. 


\bibitem[4]{ah1}  H. Ando, T. Horiuchi, E. Nakai, Weighted Hardy inequalities with infinitely many
sharp missing terms,  {\it Mathematical Journal of Ibaraki University,}  {\bf Vol. 46},  (2014), pp. 9-30. 

\bibitem[5]{BM} H. Brezis, M. Marcus, Hardy's  inequakities revisited, {\it Annali della Scuola Normale Superiore di Pisa, Classe di Scienze
  $4^e$ s\' erie, tome {\bf 25}}, {No 1-2} (1997), pp. 217-237.

\bibitem[6]{D}
E.B.Davies, The Hardy constant, Quart.J.math. Oxford (2), {\bf vol. 46}, (1995) pp. 417-431.



\bibitem[7]{dha} A. Detalla, T. Horiuchi, H. Ando, Missing terms in Hardy-Sobolev inequalities,
{\it Proceedings of the Japan Academy}, {\bf Vol. 80}, Ser. A, No. 8, 2004, pp. 160-165.
\bibitem[8]{dha2} A. Detalla, T. Horiuchi, H. Ando, Missing terms in Hardy-Sobolev inequalities and
its application, {\it Far East Journal of Mathematical Sciences}, {\bf Vol. 14}, No. 3, 2004, pp.
333-359.
\bibitem[9]{aha3}  A. Detalla, T. Horiuchi, H. Ando,  Sharp remainder terms of Hardy-Sobolev
inequalities,  {\it Mathematical Journal of Ibaraki University,} {\bf  Vol.37} (2005), pp. 39-52.  

\bibitem[10]{H}  T. Horiuchi,  
 Hardy's Inequalities  with non-doubling weights and   sharp remainders, in preparation. 

\bibitem[11]{KO}
A.Kufner and B.Opic, Hardy-type inequalities, Pitman Research notes in mathematics series
{\bf vol. 219}, [London, Longman Group UK Limited, 1990].


\bibitem[12]{MMP} M. Marcus, V. J. Mizel, Y. Pinchover, 
On  the best  constant for Hardy's inequality in {$ \bf  R^n$}, Transactions of  the American Mathematical  Society, {\bf  Vol. 350}, No. 8, August (1998), pp. 3237-3255.
\bibitem[13]{MS} T. Matskewich and P. E. Sobolevskii, The best possible constant in a generalized Hardy's inequality
for convex domains in $\bf R^n$, Nonlinear Analysis, {\bf vol. 28}, (1997) pp. 1601-1610.

\bibitem[14]{Ma} V.G. Maz'ja, Sobolev spaces (2nd edition), Springer, 2011.

\bibitem[15]{YZ} Y. Shen, Z. Chen, Sharp Hardy-Sobolev inequalities with general weights and  remainder terms, 
Journal of inequalities and applications,  Volume 2009, Article ID 419845, 24 pages
doi:10.1155/2009/419845.


\end{thebibliography}
\end{document}